\documentstyle[amssymb,amsfonts,psfig]{amsart}

\def\R{{\hbox{\bf R}}}
\def\P{{\hbox{\bf P}}}
\def\Z{{\hbox{\bf Z}}}

\def\T{{\hbox{\bf T}}}
\def\Tblift{{\overline{\hbox{\bf T}}}}

\def\Tlift{{\overline{T}}}
\def\Ylift{{\overline{Y}}}

\font \roman = cmr10 at 10 true pt

\def\be#1{\begin{equation}\label{#1}}
\def\bas{\begin{align*}}
\def\eas{\end{align*}}
\def\bi{\begin{itemize}}
\def\ei{\end{itemize}}

\def\mass{{\hbox{\roman mass}}}

\def\eps{\varepsilon}
\newenvironment{proof}{\noindent {\bf Proof} }{\endprf\par}
\def\endprf{\hfill$\square$}
\def\emph#1{{\it #1}}
\def\textbf#1{{\bf #1}}
\theoremstyle{plain}
  \newtheorem{theorem}[subsection]{Theorem}

  \newtheorem{lemma}[subsection]{Lemma}
  \newtheorem{corollary}[subsection]{Corollary}
  \newtheorem{Distance conjecture}[subsection]{Distance Conjecture}
  \newtheorem{Bilinear distance conjecture}[subsection]{Bilinear Distance Conjecture}
  \newtheorem{Furstenburg problem}[subsection]{Furstenburg problem}
  \newtheorem{Discretized Furstenburg conjecture}[subsection]{Discretized Furstenburg Conjecture}
  \newtheorem{Ring problem}[subsection]{Ring Problem}
  \newtheorem{Ring conjecture}[subsection]{Ring Conjecture}
  \newtheorem{Main theorem}[subsection]{Main Theorem}
  \newtheorem{Cauchy-Schwartz}[subsection]{Cauchy-Schwartz}

\theoremstyle{remark}

\theoremstyle{definition}
  \newtheorem{definition}[subsection]{Definition}

\include{psfig}

\begin{document}

\title[Kakeya]{New bounds for Kakeya problems}

\author{Nets Hawk Katz}
\address{Department of Mathematics, Washington University St. Louis 63130}
\email{nets@@math.wustl.edu}

\author{Terence Tao}
\address{Department of Mathematics, UCLA, Los Angeles CA 90095-1555}
\email{tao@@math.ucla.edu}


\begin{abstract} We establish new estimates on the Minkowski and Hausdorff
dimensions of Kakeya sets and we obtain new bounds on the Kakeya maximal
operator.
\end{abstract}

\maketitle

\section{Introduction}

There are many ways of formulating what is now known as the Kakeya problem.
The simplest is as follows. We define a \emph{Besicovitch set} $E \subset {\Bbb R}^n$ for $n > 1$ to be a set which contains a unit line segment in every 
direction. It is conjectured, e.g. \cite{wolff:kakeya}, \cite{borg:kakeya}, \cite{borg:high-dim}, that such a
set must have Hausdorff dimension $n$. A weaker version of the conjecture
asserts that these sets must have upper Minkowski dimension $n$.
A stronger formulation (see \cite{cordoba}, \cite{wolff:kakeya})
says that the \emph{Kakeya maximal function} 
${\cal K}_{\delta} f(\omega) := \sup_{T // \omega} {1 \over |T|} 
\int_T |f|$,
where $\delta>0$ is a fixed small number, $\omega \in S^{n-1}$, and $T$ ranges over all $1 \times \delta \dots \delta$ tubes parallel to $\omega$, is bounded
on $L^n(\R^n)$ with a bound of $C_\eps \delta^{-\eps}$ for any $\eps > 0$.  For an exposition of these problems and their applications, see \cite{Bo}, \cite{wolff:survey}, \cite{tao:notices}.

In this paper we obtain new results on all three problems in high dimensions, in the spirit of \cite{borg:high-dim} and especially \cite{katztao:+-}.  More precisely, we have

\begin{theorem} \label{summresults} 
Any Kakeya set in ${\Bbb R}^n$ must have Minkowski dimension at least
${n \over \alpha} +{\alpha-1 \over \alpha}$,
where $\alpha$ is between 1 and 2 and satisfies\footnote{Specifically, $\alpha = 1.67513\ldots$.} $\alpha^3 - 4\alpha +2=0$.
Any Kakeya set in ${\Bbb R}^n$ must have Hausdorff dimension at least
$(2-\sqrt{2})(n-4) + 3$. The Kakeya maximal function satisfies the bound
$\| {\cal K}_{\delta} f\|_{n+{3 \over 4}} \leq C_\eps \delta^{-\eps} 
 (\delta^{1-n})^{{3 \over 4n+3}} \|f\|_{{4n+3 \over 7}}$
for any $\eps > 0$ and $0 < \delta \ll 1$.
\end{theorem}

The Minkowski result is new for $n \geq 7$, improving upon the bound of $(n+2)/2 + \eps_n$ in \cite{laba:medium}; this result was previously obtained for $\alpha = 7/4$ in \cite{katztao:+-} and $\alpha = 25/13$ in \cite{borg:high-dim}. 
The Hausdorff result is new for $n \geq 5$, improving upon the bound of $(n+2)/2$ in \cite{wolff:kakeya}, and was previously obtained for $(6n+5)/11$ in \cite{katztao:+-} and $(13n+12)/25$ in \cite{borg:high-dim}.  The Hausdorff result is also superior to the stated Minkowski result for dimensions $n \leq 23$.
The maximal function result is new for $n \geq 9$, improving the bounds in \cite{wolff:kakeya} and \cite{borg:high-dim}. It is a sharp $(p,q)$ and brings $p$ to the exponent ${4n+3 \over 7}$, matching the Minkowski results of \cite{katztao:+-}. 

This paper is organized as follows.  In section 3, we present a revisionist view of  \cite{katztao:+-}, namely
Theorem \ref{conviviality}. This Theorem greatly relaxes the hypotheses
under which sums-differences lemmas can be applied, which shall be crucial in obtaining maximal function estimates (without resorting to the more involved arguments in \cite{borg:high-dim}). We also develop an iteration which gives an
easy improvement on the sums-differences lemma. This is not the best
result we have, but it serves as the model for the Hausdorff dimension
result. In section 4, we produce a more sophisticated iteration which implies
the advertised Minkowski dimension result. In section 5, we prove the maximal
estimate, and in section 6, we remove the slices in the basic iteration result
to obtain the Hausdorff bound.

We do not believe any part of Theorem \ref{summresults}
is sharp, nor that the techniques listed here are definitive.  Moreover, since we believe that the Minkowski, Hausdorff and Maximal function problems should all have the same answer, some of the ways to progress are pretty clearly indicated. We think it would greatly benefit the field if others were to take up the challenge. That means you: gentle reader!

\section{Notation}\label{notation-sec}

We write $A \gtrsim B$ if there
is a constant $C$ so that $A \geq CB$. The constant $C$ must be universal
but may vary from line to line.

We define a \emph{slope} to be an element of $\R \cup \{\infty\}$, 
and call a slope $r$ \emph{proper} if $r \neq -1$.

If $X$ is a finite set, we use $\#(X)$ to denote the cardinality of $X$.  We say that $X'$ is a \emph{refinement} of $X$ if $X' \subseteq X$ and $\#(X') \sim \#(X)$.
If $f: X \to Y$ is a map, we write $x \sim_f x'$ for 
$f(x) = f(x')$.  This induces the equivalence classes 
$[x]^{(X)}_f := \{ x' \in X: x \sim_f x' \}$.  We also define $[x]^{(X)}_{f,g} := [x]^{(X)}_{f} \cap [x]^{(X)}_g$, etc.  From Cauchy-Schwarz we record the estimate
\be{cauchy}
\#(\{(x_1,x_2) \in X \times X: s_1 \sim_f s_2\}) \geq \#(X)^2 / \#(Y).
\end{equation}

If $f_i: X \to Y_i$ are maps for $i=1, \ldots, k$ and $F: X \to Y$ is another map, we say that $F$ is \emph{determined by $f_1, \ldots, f_k$ on $X$} we have
$[x]^{(X)}_{f_1, \ldots, f_k} \subseteq [x]^{(X)}_F$
for all $x \in X$. 
If the identity map on $X$ is determined by $f_1, \ldots, f_k$, we shall say 
that $X$ is \emph{parameterized by $f_1, \ldots, f_k$}.

In practice, we will show $F$ is determined by $f_1, \ldots, f_n$ either by writing $F$ as a linear combination of the $f_i$, or by using the injectivity assumption \eqref{one-to-one}, or some combination of both.  These relationships will lead to non-trivial upper and lower cardinality bounds on various sets.

Given any sets $X_1, \ldots, X_k$ and $1 \leq i \leq k$, we define the co-ordinate functions $\gamma_i: X_1 \times \ldots \times X_k \to X_i$ by
$\gamma_i(x_1, \ldots, x_k) := x_i$.  We also define
$\gamma_{i,j}: X_1 \times \ldots \times X_k \to X_i \times X_j$ by
$\gamma_{i,j}(x_1, \ldots, x_k) := (x_i, x_j).$

We set up some notation for the ubiquitous  ``iterated popularity'' argument.  If $f: X \to Y$ is a map from one finite set $X$ to another $Y$, we define
$$ X^{<f>} := \{x \in X: \#([x]_f) \geq \#(X) / (2\#(Y)) \}.$$
We clearly have $\# (X \backslash X^{<f>}) < \#(X)/2$, and hence that $X^{<f>}$ is a refinement of $X$.  We also write $X^{<f>, <g>}$ for $(X^{<f>})^{<g>}$, etc. 

If $f(x)$ is a set-valued function on a set $X$, we write $\bigcup f(X)$ for $\bigcup_{x \in X} f(x)$.

\section{Basic iteration}\label{basic-sec}

Throughout this section we fix $Z$ to be a real vector space.  For any slope $r \in \R$, we define the maps $\pi_r: Z \times Z \to Z$ by
$\pi_r(a,b) := a+rb$ if $r \neq \infty$ and
$\pi_\infty(a,b) := b$ otherwise.  For any two slopes $r \neq r'$ we make the fundamental observation 
\be{co-ord}
Z \times Z \hbox{ is parameterized by } \pi_r, \pi_{r'}.
\end{equation}

For any two slopes $r$, $r'$, we define the double projections $\pi_{r \otimes r'}: (Z \times Z) \times (Z \times Z) \to Z \times Z$ defined 
by $\pi_{r \otimes r'}(g,g') := (\pi_r(g), \pi_{r'}(g'))$.

\begin{definition}\label{sd-def}  Let $R$ be a finite collection of proper slopes (i.e. slopes distinct from -1), and $\alpha \in \R$.  We say that the statement $SD(R, \alpha)$ holds if one has the bound
$\#(G) \lesssim \sup_{r \in R} \#(\pi_r(G))^\alpha$
whenever $G \subseteq Z \times Z$ is a finite set obeying 
\be{one-to-one}
G \hbox{ is parameterized by } \pi_{-1}
\end{equation}
(i.e. $\pi_{-1}$ is one-to-one on $G$).   
We say that the statement $SD(\alpha)$ holds if for every $\eps > 0$ there exists a finite set $R$ of proper slopes such that $SD(R, \alpha+\eps)$ holds.
\end{definition}

We remark that the $\lesssim$ in the above Definition can be automatically sharpened to $\leq$ by the iteration arguments in \cite{ruzsa} (see also \cite{katztao:+-}).
From the arguments in \cite{borg:high-dim} we see that the statement
$SD(\alpha)$ implies that Besicovitch sets have upper Minkowski dimension at least ${n \over \alpha} +{\alpha-1 \over \alpha}$; the $n$ sets $\pi_{r_i}(G)$ correspond to $n$ separate slices of the Besicovitch set.  The trivial bound $SD(\{r,r'\},2)$ for any two distinct $r,r'$ (from \eqref{co-ord}) thus implies the lower bound of $(n+1)/2$ (due to \cite{drury:xray}), whereas a result of the form $SD(1)$ would settle the Kakeya conjecture for the upper Minkowski dimension.  It is thus of interest to make $\alpha$ as low as possible, and in particular it seems not too outrageous to tentatively conjecture that $SD(1)$ is true.

In \cite{borg:high-dim} the estimate $SD(\{0,1,\infty\},2 - \frac{1}{13})$ was proven.  In \cite{katztao:+-} this was improved to $SD(\{0,1,\infty\},2 - \frac{1}{6})$, and further to

\begin{theorem}\label{012inf}\cite{katztao:+-} We have $SD(\{0,1,2,\infty\}, 2 -\frac{1}{4})$.
\end{theorem}

\begin{proof} We repeat the argument from \cite{katztao:+-}, but in a more flexible formulation. 
Define  $N := \sup_{r \in \{0,1,2,\infty\}} (\#(\pi_r(G)))$; our task is to show that $\#(G) \lesssim N^{7/4}$.

Define the set $V \subset G \times G$ of vertical line segments by 
$$V := \{ (g,g') \in G \times G: g \sim_{\pi_0} g' \} = \{( (a,b_1),(a,b_2) ): (a,b_1),(a,b_2) \in G \}.$$
We introduce the function $\nu: V \to Z$ defined by
$\nu( (a,b_1), (a,b_2) ) := a + 2b_1 - b_2$.
Since
$
\nu = \pi_2 \circ \gamma_1 - \pi_\infty \circ \gamma_2 = 2 \pi_1 \circ \gamma_1 - \pi_1 \circ \gamma_2$
we see that $\nu$ is determined by $\pi_{2 \otimes \infty}$ and also by $\pi_{1 \otimes 1}$.  From
$\nu = 2\pi_\infty \circ \gamma_1 + \pi_{-1} \circ \gamma_2$ and \eqref{one-to-one} we also see that $V$ is parameterized by $\nu, \pi_\infty \circ \gamma_1$. 

To exploit these observations we apply the iterated popularity argument.  Pick any $v_2$ in the refinement $V^{<\pi_{1 \otimes 1}>, <\pi_{2 \otimes \infty}>}$ of $V$.  By construction we have
$$ \#(\{ (v_1,v_0) \in V^{<\pi_{1 \otimes 1}>} \times V: v_2 \sim_{\pi_{2 \otimes \infty}} v_1 \sim_{\pi_{1 \otimes 1}} v_0 \})
\gtrsim \#(V)^2 / N^4.$$
For each $v_0$ there is at most one $v_1$ which contributes, by \eqref{co-ord}. 
Since $\nu$ is determined both by $\pi_{1 \otimes 1}$ and by $\pi_{2 \otimes \infty}$, we thus have $\# ([v_0]^{(V)}_\nu) \gtrsim \#(V)^2 / N^4$.
On the other hand, since $V$ is parameterized by $\nu$ and $\pi_\infty \circ \gamma_1$, we have
$\# ([v_0]^{(V)}_\nu) \leq \#(\pi_\infty \circ \gamma_1(V)) \leq N$.
Combining these two estimates we obtain
$\#(V) \lesssim N^{3/2}$.  On the other hand, from \eqref{cauchy} we have
$\#(V) \geq \#(G)^2/N$, and the claim follows.
\end{proof}

We now give a version of Theorem \ref{012inf} in which the slopes $r_i$ are more general.  Fix $G$ and a proper slope $r_0$, and let
$V = V^{r_0} := \{ (g,g') \in G^2: g \sim_{\pi_{r_0}} g' \}$.
Fix another proper slope $r_\infty \neq r_0$ and a real $s \neq 0$, and define the function $\nu: G \to \Z$ by
$$ \nu(g,g') := s \pi_{r_\infty}(g) + \pi_{-1}(g').$$
For any proper slope $r \neq r_0, r_\infty$, let $r'$ be the unique slope such that $\nu$ is determined by $\pi_{r \otimes r'}$ on $V$, or equivalently that there is an identity of the form
$$ s \pi_{r_\infty}(g) + \pi_{-1}(g') = x \pi_r(g) + y \pi_{r'}(g') + z(\pi_{r_0}(g) - \pi_{r_0}(g'))$$
for some $x,y,z \in \R$.  We refer to $r'$ as the \emph{dual slope} of $r$ with respect to $r_0$ and $\nu$.  In the special case $r_0 = 0$, $r_\infty = \infty$, $r'$ can be defined by the formula $\frac{s}{r} - \frac{1}{r'} = 1$.

Since $V$ is parameterized by $\nu$ and $\pi_{r_\infty} \circ \gamma_1$, we again have 
\be{v-upper}
\# ([v]^{(V)}_\nu) \leq \# (\pi_{r_\infty}(G)) \hbox{ for all } v \in V.
\end{equation}

On the other hand, if $r_1, r_2$ are such that the six slopes $r_1, r_2, r'_1, r'_2, -1, r_0, r_\infty$ are all distinct, then by \eqref{co-ord} $V$ is parameterized by $\pi_{r_1 \otimes r'_1}$ and $\pi_{r_2 \otimes r'_2}$.  We can thus repeat the argument in the proof of Theorem \ref{012inf} to obtain

\begin{figure}[htbp] \centering
 \ \psfig{figure=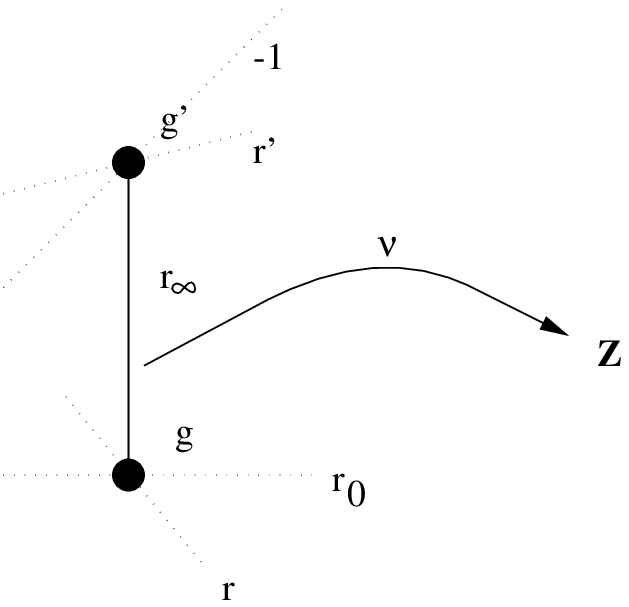,height=2in,width=2in} \psfig{figure=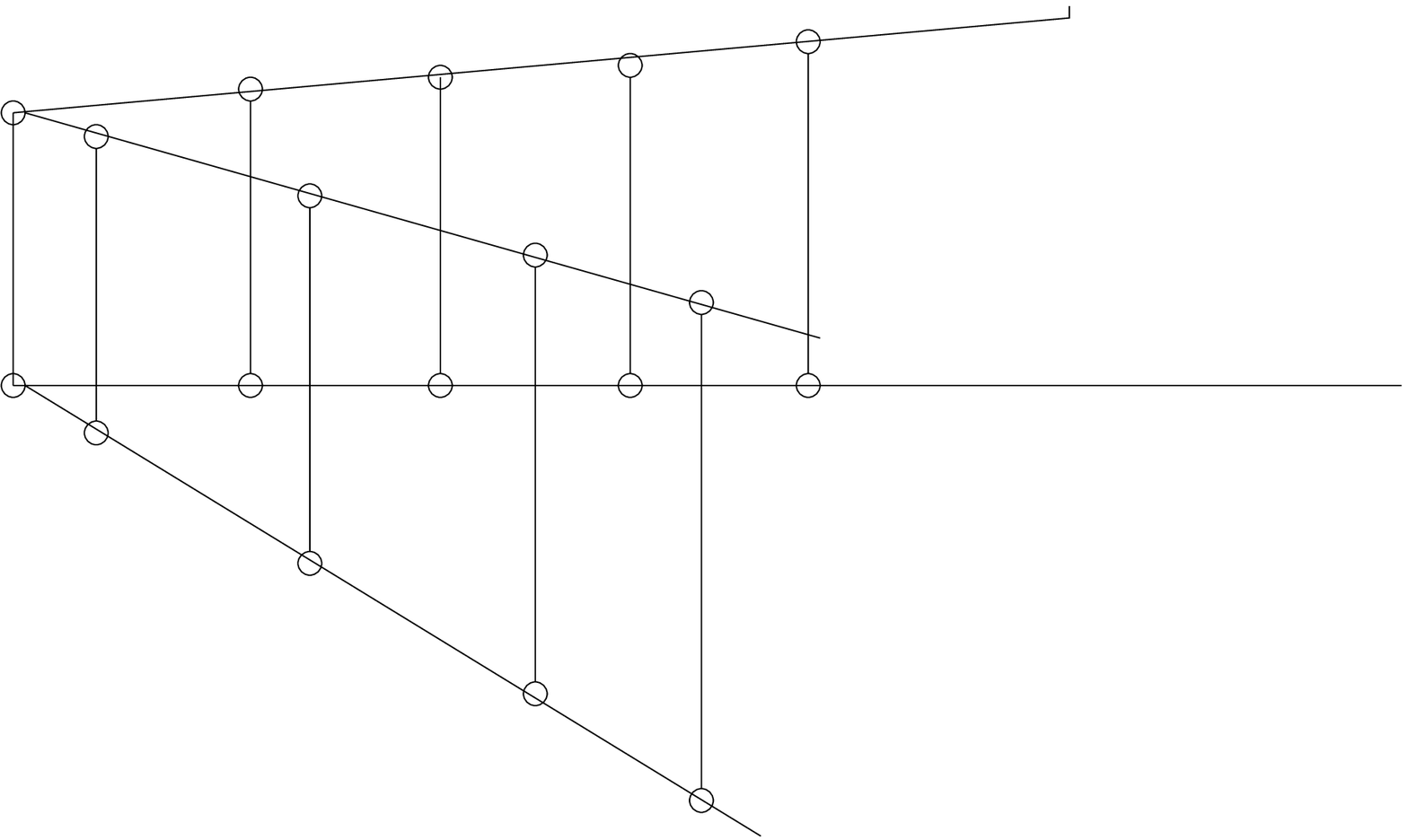,height=2in,width=2in}
 \caption{The left figure displays a vertical line segment and some of the projections $\pi_r$ one applies to this segment.   The right figure shows why we expect many vertical line segments to have the same value of $\nu$.  }
 \end{figure}

\begin{theorem}\label{conviviality} With the above assumptions, we have
$SD(\{r_0,r_1,r_1^{\prime},r_2,r_2^{\prime},r_\infty\}, {7 \over 4}).$
\end{theorem}

The disadvantage compared to Theorem \ref{012inf} is that we now need six
slices instead of four. The advantage is that the choice of slices
is much more flexible. 

We will now obtain improvements of the exponent ${7 \over 4}$, by
iterating, that is by applying sums-differences results to obtain a better lower bound for $\#([v]^{(V)}_\nu$).
The key lemma is the following:

\begin{lemma}\label{substructure} Let $G, r_0, r_\infty, V, \nu$ be as above. 
Let $r_1, \ldots, r_k$ be slopes such that the $2k+2$ slopes
\be{slope-set}
\{ r_0, r_\infty, r_1, \ldots, r_k, r'_1, \ldots, r'_k \}
\end{equation}
are all proper and disjoint. Let $N \geq 1$ be such that
$\#(\pi_r(G)) \leq N$
for all $r$ in \eqref{slope-set}.  For every $\nu_0 \in Z$, define the set $G_{\nu_0} \subseteq G$ by
$$ G_{\nu_0} := 
\gamma_1( \nu^{-1}(\nu_0) ) = \{ g \in G: (g,g') \in V \hbox{ and } \nu(g,g') = \nu_0 \hbox{ for some } g' \in G \}.$$
Then one has
\be{g-bound}
\#(G_{\nu_0}) \leq N
\end{equation}
for all $\nu_0 \in Z$.  Furthermore, there exists a $\nu_0 \in Z$ and a refinement $G'_{\nu_0}$ of $G_{\nu_0}$ such that
\be{g-smallproj}
\pi_{r_j}(G'_{\nu_0}) \lesssim \frac{\#(G_{\nu_0})}{\#(V)/N^2}
\end{equation}
for $j = 1, \ldots, k$.  (The implicit constants depend on $k$).
\end{lemma}

\begin{proof} 
Since $\nu$ is determined by $\pi_{r_\infty} \circ \gamma_1$ on $V$, it is determined by $\gamma_1$, and we thus have
$\#(G_{\nu_0}) = \# (\nu^{-1}(\nu_0) ).$
The claim \eqref{g-bound} then follows from \eqref{v-upper}.

Define 
$V' := V^{<\pi_{r_1 \otimes r'_1}>, \ldots, <\pi_{r_k \otimes r'_k}>}$.
Since $V'$ is a refinement of $V$, we may find a $\nu_0 \in \Z$ such that
$\#(V' \cap \nu^{-1}(\nu_0)) \sim \#(\nu^{-1}(\nu_0)) \neq 0.$
Fix this $\nu_0$, and define $G'_{\nu_0} \subseteq G_{\nu_0}$ by
$G'_{\nu_0} := \{ g \in G: (g,g') \in V' \cap \nu^{-1}(\nu_0) \hbox{ for some } g' \in G \}$.
Since $\nu$ is determined by $\gamma_1$, we see that
$$ \# (G'_{\nu_0}) = \# (V' \cap \nu^{-1}(\nu_0)) \sim \#(\nu^{-1}(\nu_0))$$
so $G'_{\nu_0}$ is a refinement of $G_{\nu_0}$.

Now we prove \eqref{g-smallproj}.  Fix $j$.  Since $\nu$ is determined by $\pi_{r,r'}$ on $V$, we have some identity of the form
$\nu(g,g') = x \pi_r(g) + y \pi_{r'}(g')$
for some non-zero scalars $x,y$ and all $(g,g') \in V$.  In particular, we see that $\pi_{r,r'}$ is determined by $\pi_r \circ \gamma_1$ on $V^0_{\nu_0}$.  Thus to prove \eqref{g-smallproj} it suffices to show that
$$ \pi_{r_j \otimes r'_j}(V' \cap \nu^{-1}(\nu_0)) \lesssim \frac{\#(\nu^{-1}(\nu_0))}{\#(V)/N^2}.$$
On the other hand, from the construction of $V'$ we have
$$ \#(\{ v' \in \nu^{-1}_{\nu_0}: v \sim_{\pi_{r_j \otimes r'_j}} v' \}) 
\gtrsim \#(V)/N^2$$
for all $v \in V' \cap \nu^{-1}(\nu_0)$.  The claim follows. 
\end{proof}

\begin{corollary}\label{iteration}  For any $1 < \beta \leq 2$, we have
$SD(\beta) \implies SD({4\beta -1 \over 2\beta}).$
\end{corollary}

\begin{proof} By an obvious limiting argument it suffices
to show that for any finite set of proper slopes $R_0$, we can find a finite set $R$ of proper slopes such that
$$SD(R_0,\beta) \implies SD(R, {4\beta -1 \over 2\beta}).$$
Fix $R_0$.  We now pick $r_0$, $r_\infty$, and $\nu$ in such a way
that the slopes in \eqref{slope-set} are proper and distinct. It is clear that this distinctness property holds for generic choices of $r_0$, $r_\infty$, and $\nu$.  We then set
$R$ to equal to \eqref{slope-set}.
Let $N \geq 1$, and let $G \subseteq Z \times Z$ obey \eqref{one-to-one} be such that $\#(\pi_r(G)) \leq N$ for all $r \in R$.  Then we may apply Lemma \ref{substructure} to obtain a $\nu_0 \in Z$ and $G'_{\nu_0} \subseteq G$ obeying the conclusions of that Lemma. 

We now apply the hypothesis $SD(R_0,\beta)$ with $G$ replaced by the smaller set $G'_{\nu_0}$ (note that \eqref{one-to-one} is still true).  From \eqref{g-smallproj} we have
$
\# (G_{\nu_0}) \lesssim (\frac{\#(G_{\nu_0})}{\#(V)/N^2})^\beta$;
combining this with \eqref{g-bound} we obtain after some algebra
$\#(V) \lesssim N^{{3 \beta-1 \over \beta}}$.
Combining this from the bound $\#(V) \geq \#(G)^2/N$ from \eqref{cauchy} we obtain the claim. 
\end{proof}

The reader should observe that proof of Theorem \ref{012inf} is essentially a special case of the proof of Corollary \ref{iteration}, specialized to $\beta = 2$.  If we instead iterate Corollary \ref{iteration} and solving for the fixed point of $\beta = \frac{4\beta-1}{2\beta}$ we obtain
$SD(1+{\sqrt{2} \over 2})$.  This implies a bound on the Minkowski problem, but this will be superceded by a more sophisticated iteration argument in the next section.  However, we shall use the ideas used to prove this result in Section \ref{hausdorff-sec}.

\section{Advanced Iteration}

The goal of this section is to prove $SD(\alpha)$, where $\alpha$ is as in 
Theorem \ref{summresults}.

Just as $SD(1 + \frac{\sqrt{2}}{2})$ follows easily from iterating Corollary \ref{iteration}, $SD(\alpha)$ follows from iterating

\begin{theorem}\label{was27-thm}  For any $1 < \beta \leq 2$, we have
$
 SD(\beta) \implies SD({3\beta^2+2\beta-2 \over 
\beta^2+3\beta-2})
$.
\end{theorem}

Before we give the rigorous proof of this Theorem, we first give a

\emph{Heuristic Proof.} In order to prove the Theorem, we must study more sophisticated objects than vertical line segments, namely corners.  Let $r_1, r_2$ be proper slopes, and define the set $C = C^{r_1,r_2}$ of \emph{corners with slopes $r_1, r_2$} by
\be{c-def}
C^{r_1,r_2} := \{(g_1,g_2,g_3) \in G: g_1 \sim_{\pi_{r_1}} g_2 \sim_{\pi_{r_2}} g_3 \}.
\end{equation}

\begin{figure}[htbp] \centering
 \ \psfig{figure=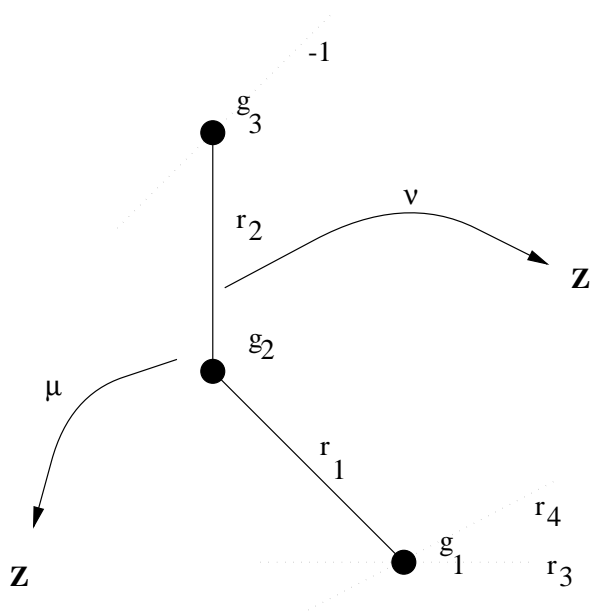}
 \caption{A corner, and some of the projections $\pi_r$ one
applies to this corner.   The map $\mu$ takes corners to
elements of $Z$.}
 \end{figure}

On $C$ we have the projections $\gamma_i: C \to G$ for $i=1,2,3$, as well as the projections $\gamma_{1,2}: C \to V^{(r_1)}$ and $\gamma_{2,3}: C \to V^{(r_2)}$.  We shall also need proper slope $r_3$ and a map $\mu =  \mu^{(r_3)}: C \to Z$, defined by $\mu(g_1,g_2,g_3) := \pi_{r_3}(g_1) + \pi_{-1}(g_3)$.
In this argument $C$ and $\mu$ shall roughly play the role of $V$ and $\nu$ in the previous section.

Assume that $\#(\pi_s(G)) \leq N$ for all proper slopes $s$ which we shall use.
If $V^{(r)}$ be the set of vertical line segments with slope $r$, we shall make the heuristic uniformity assumption
\be{2pi}
\#([v]_{\pi_{s_1 \otimes s_2}} )\lesssim 
\#(V^{(r)}) / N^2 \sim \#(G)^2/N^3
\end{equation}
for any $v \in V^{(r)}$ and any $s_1, s_2$ we shall need.
Moreover, we shall make the assumption that for slope $r$ and any choice of
function $\nu$ as in the previous section, we assume that the equivalence classes $[v]_\nu$ in $V^{(r)}$ for $v \in V^{(r)}$ all have roughly the same cardinality, which we will refer by abuse of notation as
$\#([v]_\nu)$.  In particular, we assume
\be{nu-bound}
\#(\nu(V^{(r)})) \lesssim \frac{\#(G)^2/N}{\#([v]_\nu)}
\end{equation}
for all $\nu, r$.
It is largely because of these assumptions that the current proof
is merely heuristic. 

From the definition of $\mu$ we see that
$\pi_{-1} \circ \gamma_3$ is determined by $\mu, \pi_{r_3} \circ \gamma_1$
 on $C$, 
so by \eqref{one-to-one} we have
\be{l-3}
\gamma_3 \hbox{ is determined by } \mu, \pi_{r_3} \circ \gamma_1 \hbox{ on C}.
\end{equation}
On the other hand, from \eqref{c-def}, \eqref{co-ord} we see that
$\gamma_2$ is determined by $\gamma_1, \gamma_3$ on $C$
so that $C$ is parameterized by $\gamma_1, \mu$.
In particular, for any $c \in C$ we have
$\#([c]_\mu) = \#(\gamma_1([c]_\mu))$.
We shall now obtain upper and lower bounds for this quantity for generic values of $c$.

To obtain the upper bound, we first consider the equivalence class $[c]_{\pi_{r_3} \circ \gamma_1, \mu}$.  From \eqref{l-3} we have
$
\#([c]_{\pi_{r_3} \circ \gamma_1, \mu}) = \#(\gamma_{1,2}([c]_{\pi_{r_3} \circ \gamma_1, \mu}))$. 
From \eqref{l-3} and \eqref{c-def} we see that $\pi_{r_3 \otimes r_2} \circ \gamma_{1,2}$ is determined by $\pi_{r_3} \circ \gamma_1$ and $\mu$ on $C$. From \eqref{2pi} we thus have
$\# ([c]_{\pi_{r_3} \circ \gamma_1, \mu}) \lesssim \#(G)^2/N^3$.
Since there are at most $N$ possible values of $\pi_{r_3} \circ \gamma_1$, we have shown that
\be{was21} \#(\gamma_1([c]_\mu)) = \# [c]_\mu \lesssim \#(G)^2/N^2. \end{equation}

To obtain the lower bound on $\#(\gamma_1([c]_\mu))$, we shall apply the hypothesis $SD(\beta)$. To do this, we must bound the cardinality of $\pi_{r_4} \circ \gamma_1([c]_\mu)$ for a typical $r_4$. The main point is that for a fixed (generic) $r_4$, there exists a function $\nu: V^{(r_2)} \to Z$ of the type discussed in the previous section,
such that $\mu$, $\pi_{r_4} \circ \gamma_1$, and $\nu \circ \gamma_{2,3}$ are linearly dependent, and in particular that
\be{nu-mu}
[c]_{\pi_{r_4} \circ \gamma_1, \mu} = [c]_{\pi_{r_4} \circ \gamma_1, \nu \circ \gamma_{2,3}}
\end{equation}
for all $c \in C$.
Since $\pi_r \circ \gamma_1$ maps to a set of cardinality at most $N$, and $\nu \circ \gamma_{2,3}$ maps to a set of cardinality at most $\frac{\#(G)^2/N}{\#([v]_\nu)}$ by \eqref{nu-bound}, we thus expect 
$$ \#([c]_{\pi_{r_4} \circ \gamma_1, \mu}) \gtrsim \frac{\#(C)}{N ((\#(G)^2/N)/\#([v]_\nu))}.$$ 
From the lower bound 
$ \#(C) \gtrsim \frac{\#(G)^3}{N^2}$ (cf. \cite{katztao:+-}, Lemma 2.1) we thus expect
$$\#([c]_{\pi_{r_4} \circ \gamma_1, \mu}) \gtrsim \#(G) \#([v]_\nu) / N^2, $$
and hence that
$$ \#(\pi_{r_4} \circ \gamma_1([c]_\mu)) \lesssim \frac{\#([c]_\mu)}{ \#(G) \#([v]_\nu) / N^2 }.$$
We apply $SD(\beta)$ to this estimate and \eqref{was21} to obtain (ignoring epsilons)
$$
({\#(G)\#([v]_\nu) \over  N^2})^{{\beta \over \beta-1}}
\lesssim {\#(G)^2 \over N^2}. $$
However, from the iteration scheme of the previous section, we
have
\be{was25}
\#([v]_\nu) \gtrsim ( {\#(V) \over N^2})^{{\beta \over \beta -1}}
\gtrsim ( {\#(G)^2 \over N^3})^{{\beta \over \beta -1}}. 
\end{equation}
Combining these estimates we obtain 
$\#(G) \lesssim N^{{3\beta^2+2\beta-2 \over \beta^2+3\beta-2}}.
$ as desired.
\endprf

We now discuss what is needed to make this heuristic proof rigorous.  Most of the assumptions can be formalized rigorously using standard uniformization tricks such as the iterated popularity arguments of the previous section. However, a difficulty arises when trying to derive \eqref{2pi}.  If one has breakdown of uniformity here then we suffer a loss in \eqref{was21}.  In principle we can recover this because we will also gain in \eqref{was25}, except that the slopes used in \eqref{was25} are not the same as those for \eqref{was21}. We will resolve this issue by being prepared to swap the choices of slices $r$ a finite number of times (losing an epsilon in the exponents each time); cf. the uniformization arguments in \cite{wolff:xray}.

We begin the rigorous argument.  By standard limiting arguments it suffices to show

\begin{theorem}  Let $1 < \beta \leq 2$, and let $R_0$ be a finite set of proper slopes such that $SD(R_0, \beta)$ holds.  Then for any integer $M \gg 1$ we can find a finite set $R$ of proper slopes and $C(R_0) > 0$ such that 
$SD(R, {3\beta^2+2\beta-2 \over \beta^2+3\beta-2} + \frac{C(R_0)}{M})$
holds.
\end{theorem}

\begin{proof}
Fix $\beta$, $R_0$, $M$.  The first step is to construct the set $R$ of slopes.  The construction will be quite involved, but is necessary in order to perform the rest of the argument rigorously.

Let $r_1, r_2, r_3$ be distinct proper slopes.  We can then form as before the set $C^{r_1,r_2}$ of corners, and the map $\mu^{(r_3)}$ on $C^{r_1,r_2}$.  Given any other slope $r_4$, recall that there exists a function $\nu = \nu_{r_4}$ on $V^{(r_2)}$ such that \eqref{nu-mu} holds.  For all but a finite number of exceptional $r_4$, this $\nu$ is of the type discussed in the previous section, and in particular for any fixed such $r_4$, and all but a finite number of proper slopes $r$, there exists a proper slope $r'$ (distinct from $r_1, r_2, r_3, r_4, r$) such that $\nu$ is a linear combination of $\pi_r \circ \gamma_1$ and $\pi_{r'} \circ \gamma_2$.

It is too unreasonable to expect all of these slopes to lie in $R_0$.  On the other hand, from Definition \ref{sd-def} and the hypothesis $SD(R_0, \beta)$ we have $SD(L(R_0), \beta)$ whenever $L$ is a fractional linear transformation on $\R$ which preserves -1.  This gives us much more flexibility since we can always choose such an $L$ even if a finite number of exceptional $L$ are somehow prohibited.  For instance, for fixed $r_1, r_2, r_3$, and for $r_4$ avoiding a finite number of exceptional values (depending on $r_1, r_2, r_3$), one can find a set $R_{r_1,r_2,r_3,r_4}$ of proper slopes such that $SD(R_{r_1,r_2,r_3,r_4},\beta)$ holds, and such that a dual proper slope $r'$ distinct from $r_1, r_2, r_3, r_4, r$ exists for each $r \in R_{r_1,r_2,r_3,r_4}$.  Similarly, for fixed $r_1, r_2, r_3$, there exists a set $R_{r_1,r_2,r_3}$ of proper slopes $r_4$ avoiding $r_1$, $r_2$, $r_3$, and the exceptional values mentioned earlier, such that $SD(R_{r_1,r_2,r_3},\beta)$ holds.

Fix the sets $R_{r_1, r_2,r_3}$ and $R_{r_1,r_2,r_3,r_4}$.  For any three distinct proper slopes $r_1, r_2, r_3$, we define the set of slopes
$$R^*(r_1,r_2,r_3) := \{r_1,r_2,r_3\} \cup R_{r_1,r_2,r_3} \cup
\bigcup_{r_4 \in R_{r_1,r_2,r_3}} \bigcup_{r \in R_{r_1,r_2,r_3,r_4}} \{r,r'\}$$
and set $\T(r_1,r_2,r_3)$ of triples of distinct proper slopes by
$$ \T(r_1,r_2,r_3) := \{(r_2, r, r'): r_4 \in R_{r_1,r_2,r_3}, r \in R_{r_1,r_2,r_3,r_4} \}.$$

Let $T_0 := \{ (0, 1, 2) \}$ (for instance), and define recursively $T_1, T_2, \ldots, T_M$ by $T_{j+1} := \bigcup \T(T_j)$
for $j = 1, \ldots, M$; note that the $T_j$ always consist of triples of three distinct proper slopes.  We then set $R$ equal to
$R := \bigcup_{j=0}^{M} \bigcup R^*(T_j)$.

The idea will be to try to run the heuristic argument using some triples  from the set $T_M$.  If at least one of these triples satisfies a certain uniformity property then the argument will run smoothly.  If all the triples from $T_M$ fail to be uniform, then we pass to the triples $T_{M-1}$.  The point is that the failure of uniformity for $T_M$ will allow us to be less strict about the uniformity required for $T_{M-1}$ (we can lose an additional factor of $N^{C/M}$ or so).  We repeat this process as long as necessary.  In the worst case we fall back all the way to $T_0$, but the uniformity requirement is now trivial
(cf. \cite{wolff:xray}, \cite{laba:x-ray}).

Henceforth all implicit constants in the $\lesssim$ notation will be allowed to depend on $\beta$, $M$, $R$, $R_0$, while constants denoted by $C$ are only allowed to depend on $R_0$.

Let $G$ obey \eqref{one-to-one} and set $N := \sup_{r \in R} \#(\pi_r(G))$; we may assume $N \gg 1$. Our task is to show that
\be{g-task}
\#(G) \lesssim N^{C/M} N^{{3\beta^2+2\beta-2 \over 
\beta^2+3\beta-2}}.
\end{equation}
Define refinement $G' := G^{<r_1>, <r_2>, \ldots, <r_k>}$ of $G$, where $r_1, \ldots, r_k$ is an arbitrary enumeration of the finite set $R$.  By construction we can choose for each $g \in G'$ and $r \in R$ a set $X_{g,r} \subseteq [g]^{(G)}_{\pi_r}$ such that
$\#(X_{g,r}) \sim \#(G)/N$. 

Fix the sets $X_{g,r}$.  We define the modified vertical line segment sets $\tilde V^{(r)} \subset V^{(r)}$ for $r \in R$ by
$\tilde V^{(r)} := \{ (g,g'): g \in G', g' \in X_{g,r} \}$.
Clearly we have 
\be{tvr-card}
\#(\tilde V^{(r)}) \sim \#(G)^2/N.
\end{equation}

Let $0 \leq k \leq M$, and let $(r_1,r_2,r_3)$ be an element of $T_k$.  We say that $(r_1,r_2,r_3)$ is \emph{$k$-uniform} if we have
$$ \#(\{ v \in \tilde V^{(r_1)}: \#([v]^{(\tilde V^{(r_1)})}_{\pi_{r_2 \otimes r_3}})
< \frac{\#(\tilde V^{(r_1)})}{N^2} N^{\rho_k}
\}) \geq N^{-\frac{1}{M}-\frac{k}{M^2}} \#(\tilde V^{(r_1)})$$
where $\rho_k := 100(1-k/M)$.  For $k$ close to $M$, this property asserts that the map $\pi_{r_2 \otimes r_3}$ maps $\tilde V^{(r_1)}$ maps evenly onto $\pi_{r_2}(G) \times \pi_{r_3}(G)$.  However this property becomes weaker as $k$ decreases.  For instance, it is clear that the singleton element of $T_0$ is $0$-uniform.  

With the notation of the previous paragraph, we say that $(r_1,r_2,r_3) \in T_k$ is \emph{$k$-chunky} if there exists a subset $\tilde V^{(r_1,r_2,r_3)}$ of $\tilde V^{(r_1)}$ such that
$\#(\pi_{r_2 \otimes r_3}(\tilde V^{(r_1,r_2,r_3)})) 
\lesssim N^{2 - \rho_k}$ and
$$ \# (\tilde V^{(r_1,r_2,r_3)}) \geq (1 - C N^{-\frac{1}{M}-\frac{k}{M^2}}) \#(\tilde V^{(r_1)}).$$
Clearly if a triple $(r_1,r_2,r_3) \in T_k$ fails to be $k$-uniform, then it is $k$-chunky (just set $\tilde V^{(r_1,r_2,r_3)}$ equal to the appropriate elements of $\tilde V^{(r_1)}$).  Also, every triple in $T_M$ is trivially $M$-chunky.

From the above discussion, it is clear that we can find a $0 \leq k < M$ and a triple $(r_1,r_2,r_3) \in T_k$ which is $k$-uniform, and such that every triple in $T_{k+1}$ is $k+1$-chunky.  Henceforth $k$ and $(r_1,r_2,r_3)$ are fixed.
The main issue here is to ensure the powers of $N^{\rho_k}$ one loses in the $k$-uniformity property will be compensated for by the gains in $N^{\rho_k}$ one will obtain from the $k+1$-chunkiness property.

For each $(r_2,r,r') \in \T(r_1,r_2,r_3)$, we have $(r_2,r,r') \subseteq T_{k+1}$, and hence that $(r_2,r,r')$ is $k+1$-chunky.  Hence we have
$$ \#(\tilde V^{r_2} \backslash \tilde V^{(r_2,r,r')})
\lesssim N^{-\frac{1}{M}-\frac{k+1}{M^2}} \#(\tilde V^{(r_2)})$$
and
\be{v-smallproj}
\#(\pi_{r \otimes r'}(\tilde V^{(r_2,r,r')})) 
\lesssim N^{2-\rho_k}.
\end{equation}
We introduce the refinement $V' := \bigcap_{(r_2,r,r') \in \T(r_1,r_2,r_3)} \tilde V^{(r_2,r,r')}$ of $\tilde V^{(r_2)}$ and observe
\be{jump}
\#(\tilde V^{(r_2)} \backslash V') \lesssim N^{-\frac{1}{M}-\frac{k+1}{M^2}} \#(\tilde V^{(r_2)}).
\end{equation}
We now run a modified iterated popularity argument.  For any map $f: X \to Y$ between finite sets, define 
$$ X^{<<f>>} := \{ x \in X: \#([x]_f) \geq N^{-\frac{100}{M}} \#(X) / \#(Y) \}.$$
Observe the strong refinement property
$\#(X \backslash X^{<<f>>}) \lesssim N^{-\frac{100}{M}} \#(X)$.
We now define the refinement $V'' := (V')^{<<f_1>>, \ldots, <<f_s>>}$ of $V'$,
where $f_1, \ldots, f_s$ is some arbitrary enumeration of the projections 
$\{ \pi_{r \otimes r'}: (r_2,r,r') \in \T(r_1,r_2,r_3) \},$
and the projections $\pi_{r \otimes r'}$ are thought of as mapping to $\pi_{r \otimes r'}(\tilde V^{(r_2,r,r')})$.

From the strong refinement property we have (for $N$ sufficiently large)
$\#(V' \backslash V'') \lesssim N^{-\frac{100}{M}} \#(V')$;
combining this with \eqref{jump} we obtain
\be{jump-2}
\#(\tilde V^{(r_2)} \backslash V'') \lesssim N^{-\frac{1}{M}-\frac{k+1}{M^2}} \#(\tilde V^{(r_2)}).
\end{equation}

Let $r_4$ be a slope in $R_{r_1,r_2,r_3}$, and let $\nu$ be the associated map on $\tilde V^{(r_2)}$.  Introduce the sets
$V''(r_4) := \{ v \in V'': \#([v]^{(V'')}_\nu) \sim \#([v]^{(\tilde V^{(r_2)})}_\nu) \}$.
From \eqref{jump-2} we see that
$V''(r_4) \lesssim N^{-\frac{1}{M}-\frac{k+1}{M^2}} \#(\tilde V^{(r_2)})$.
If we then introduce 
$V''' := V'' \backslash \bigcup V''(R_{r_1,r_2,r_3})$ of $V''$ then
\be{jump-3}
\#(\tilde V^{(r_2)} \backslash V''') \lesssim N^{-\frac{1}{M}-\frac{k+1}{M^2}} \#(\tilde V^{(r_2)}).
\end{equation}

We shall need a bound on the cardinality of $\nu(V''')$.  To do this we shall first repeat the argument in Corollary \ref{iteration}.  By the construction of $V''$ we have
$$ \#([v]^{(\tilde V^{(r_2)})}_{\pi_{r \otimes r'}}) \gtrsim N^{-C/M} \frac{\#(V'')} {\#(\pi_{r \otimes r'}(\tilde V^{(r_2,r,r')}))}$$
for all $v \in V'''$;
by \eqref{v-smallproj}, \eqref{jump-2}, \eqref{tvr-card} we thus have
$$ \#([v]^{(\tilde V^{(r_2)})}_{\pi_{r \otimes r'}}) \gtrsim N^{-C/M} \frac{\#(G)^2/N}{ N^{2-\rho_k} }.$$
Since $\nu$ is determined by $\pi_{r \otimes r'}$, we thus have
$$
 \pi_{r \otimes r'}([v]^{(V''')}_\nu) \lesssim N^{C/M} \frac{N^{2-\rho_k}}{\#(G)^2/N} \#([v]^{(\tilde V^{(r_2)})}_\nu).
$$
Fix $v$, and define the set $G' \subseteq G$ by
$G' := \gamma_1([v]^{(V''')}_\nu)$.
Since $[v]^{(V''')}_\nu$ is parameterized by $\gamma_1$, we see from construction that
\be{sweet}
\#([v]^{(\tilde V^{(r_2)})}_\nu) \sim  \#([v]^{(V'')}_\nu)
= \#([v]^{(V''')}_\nu) = \#(G').
\end{equation}
Inserting this into the previous we obtain
$V^{(r_2)}$, we thus see from \eqref{sweet} that
$$ \pi_r(G') \leq \pi_{r \otimes r'}([v]^{(V''')}_\nu) \lesssim N^{C/M} \frac{N^{3-\rho_k}}{\#(G)^2} \#(G').$$
Applying $SD(R_{r_1,r_2,r_3,r_4},\beta)$ we thus obtain that
$\#(G') \lesssim (N^{C/M} \frac{N^{3-\rho_k}}{\#(G)^2} \#(G'))^\beta$.
By \eqref{sweet} we thus have 
$$ \#([v]^{(V''')}_\nu) \gtrsim N^{-C/M} (\frac{N^{3-\rho_k}}{\#(G)^2})^{-\frac{\beta}{\beta-1}}.$$
Since this holds for all $v \in V'''$, we thus obtain
$$ \#(\nu(V''')) \lesssim N^{C/M} (\frac{N^{3-\rho_k}}{\#(G)^2})^{\frac{\beta}{\beta-1}} \#(V'''),$$
which simplifies using \eqref{jump-3}, \eqref{tvr-card} to
\be{vppp-nu}
\#(\nu(V''')) \lesssim N^{C/M} (\frac{N^{3-\rho_k}}{\#(G)^2})^{\frac{\beta}{\beta-1}} \frac{\#(G)^2}{N}.
\end{equation}

Having defined $V'''$, we now introduce the set
$C := \{ (g_1, g_2, g_3) \in G^3: (g_2,g_1) \in \tilde V^{(r_1)}; (g_2,g_3) \in \tilde V^{(r_2)}\}$ of corners.  Since $(r_1,r_2,r_3)$ is $k$-uniform, we see that
$$ \# (\{ v \in \tilde V^{(r_1)}: [v]^{(\tilde V^{(r_1)})}_{\pi_{r_2 \otimes r_3}} < \frac{\#(\tilde V^{(r_1)})}{N^2} N^{\rho_k} \})
\geq N^{-\frac{1}{M}-\frac{k}{M^2}} \#(\tilde V^{(r_1)}).$$
From the construction of $X_{g,r}$ we thus have
$$ \# (\{ c \in C: [\gamma_{2,1}(c)]^{(\tilde V^{(r_1)})}_{\pi_{r_2 \otimes r_3}} < \frac{\#(\tilde V^{(r_1)})}{N^2} N^{\rho_k} \})
\gtrsim N^{-\frac{1}{M}-\frac{k}{M^2}} \#(C).$$
If we therefore define the set $C'$ by
\be{cp-def}
C' := \{ c \in C: [\gamma_{2,1}(c)]^{(\tilde V^{(r_1)})}_{\pi_{r_2 \otimes r_3}} < \#(\tilde V^{(r_1)}) / N^{2-\rho_k}; \gamma_{2,3}(c) \in V''' \}
\end{equation}
then we see from \eqref{jump-3} that
$\#(C') \gtrsim N^{-\frac{1}{M}-\frac{k}{M^2}} \#(C).$
From the construction of $C$ we thus have 
\be{cp-bound}
N^{-C/M} \frac{\#(G)^3}{N^2} \lesssim \#(C') \lesssim \frac{\#(G)^3}{N^2}.
\end{equation}

For any $r_4 \in R_{r_1,r_2,r_3}$, define the map $f_{r_4}: C' \to \pi_{r_4}(G) \times Z$ by $f_{r_4}(g_1,g_2,g_3) := (\pi_1(g_1), \nu(g_2,g_3)),$
where $\nu: V^{(r_2)} \to Z$ is the map associated to $r_4$. 
We once again apply the iteration arguments of Lemma \ref{substructure},  defining the refinement $C'' := (C')^{<f_1>, <f_2>, \ldots, <f_s>}$ of $C'$,
where $f_1, \ldots, f_s$ is an arbitrary enumeration of the functions $\{ f_{r_4}: r_4 \in R_{r_1,r_2,r_3} \}$. Since $C''$ is a refinement of $C'$,
we may fix a $c \in C''$ such that
$[c]^{(C'')}_\mu$ is a refinement of $[c]^{(C')}_\mu$.
Let $G'' \subseteq G$ denote the set
$G'' := \gamma_1([c]^{(C'')}_\mu)$. We now obtain upper and lower bounds for the size of $G''$.

To obtain lower bounds, we argue as in Lemma \ref{substructure}.  From the construction of $C''$ we have
$\#([c']^{(C')}_{f_{r_4}})  \gtrsim \#(C'') / \#(f_{r_4}(C'))$
for any $c' \in [c]^{(C'')}_\mu$.  Using this and the fact that $\mu$ is determined by $f_{r_4}$ we obtain
$\#(f_{r_4}( [c]^{(C'')}_\mu )) \lesssim (\#(f_{r_4}(C'))/\#(C'')) \#([c]^{(C'')}_\mu )$.
On the other hand, from \eqref{vppp-nu} we have
$$
\#(f_{r_4}(C')) \lesssim N \times N^{C/M} (\frac{N^{3-\rho_k}}{\#(G)^2})^{\frac{\beta}{\beta-1}} \frac{\#(G)^2}{N}.
$$
Combining these bounds with \eqref{cp-bound} we obtain
$$ \#(f_{r_4}( [c]^{(C'')}_\mu )) \lesssim N^{C/M} N^{-\frac{\beta}{\beta-1}\rho_k + \frac{3\beta}{\beta-1} + 2} \#(G)^{-\frac{2}{\beta-1} - 3}
 \#([c]^{(C'')}_\mu ).$$
Since $[c]^{(C'')}_\mu$ is parameterized by $\gamma_1$, we have
\be{flip}
 \#(f_{r_4}( [c]^{(C'')}_\mu )) = \#(\pi_{r_4}( G'' )) \hbox{ and }
\#([c]^{(C'')}_\mu = \#(G'').
\end{equation}
Applying $SD(R_{r_1,r_2,r_3,r_4},\beta)$ we obtain the lower bound
\be{advanced}
\#(G'') \gtrsim N^{-C/M} \left(N^{-\frac{\beta}{\beta-1}\rho_k + \frac{3\beta}{\beta-1} + 2} \#(G)^{-\frac{2}{\beta-1} - 3}\right)^{-\frac{\beta}{\beta-1}}.
\end{equation}

Now we obtain upper bounds on $\# ([c]^{(C'')}_\mu)$. 
We observe as in the Heuristic proof that $\pi_{r_2 \otimes r_3} \circ \gamma_{2,1}$ is determined by $\mu$ and $\pi_{r_3} \circ \gamma_1$, so we have
$$ 
\# ([c]^{(C'')}_{\mu, \pi_{r_3} \circ \gamma_1})
=
\# ([c]^{(C'')}_{\mu, \pi_{r_2 \otimes r_3} \circ \gamma_{2,1} })
\leq
\# ([\gamma_{2,1}(c)]^{(\tilde V^{(r_1)})}_{\pi_{r_2 \otimes r_3} }).$$
By \eqref{cp-def}, \eqref{tvr-card} we thus have
$$ \# ([c]^{(C'')}_{\mu, \pi_{r_3} \circ \gamma_1}) \lesssim \#(\tilde V^{(r_1)}) / N^{2-\rho_k} \sim \#(G)^2 / N^{3 - \rho_k}.$$
Since the range of $\pi_{r_3} \circ \gamma_1$ has cardinality at most $N$, we thus have $\# ([c]^{(C'')}_{\mu}) \lesssim \#(G)^2 / N^{2-\rho_k}$.
Combining this with \eqref{advanced} and \eqref{flip} we obtain
$$
N^{-C/M} \left(N^{-\frac{\beta}{\beta-1}\rho_k + \frac{3\beta}{\beta-1} + 2} \#(G)^{-\frac{2}{\beta-1} - 3}\right)^{-\frac{\beta}{\beta-1}}
\lesssim
\#(G)^2 / N^{2 - \rho_k}.$$
which simplifies using $\rho_k \geq 0$ to \eqref{g-task} as desired.
\end{proof}

\section{Maximal Function}

The purpose of this section is to apply the ideas of
Theorem \ref{conviviality} to obtain
new sharp $(p,q)$ bounds on the Kakeya maximal function, defined in the introduction.

In this section $0 < \eps \ll 1$ is a small fixed parameter, and $N \gg 1$ is a large parameter.  We use $A \lessapprox B$ to 
denote the estimate $A \leq C_\eps N^{C\eps} B$ for some constants $C_\eps$, $C$, and $A \approx B$ to denote $A \lessapprox B \lessapprox A$.  We use $\delta$ to denote the quantity $\delta := 1/N$.  We redefine a \emph{refinement} $X'$ of a set $X$ to be any subset $X' \subseteq X$ such that $\#(X') \approx \#(X)$.

The objective of this section is to prove the estimate $\| {\cal K}_{\delta} f\|_{(n-1)p'} \lessapprox  \delta^{\frac{n}{p}-1} \|f\|_{p}$ when $1 \leq p \leq (4n+3)/7$; the reader may verify that the exponents in the estimate are sharp.
We will assume $n > 8$ since the claim follows from the $p=(n+2)/2$ estimate of \cite{wolff:kakeya} otherwise.  This estimate implies the weaker Hausdorff and upper Minkowski dimension in \cite{katztao:+-}, but not the other results in this paper.
 
We may take $p := (4n+3)/7$.  To prove the above estimate we first make some standard reductions.  We first observe that we may as well restrict $f$ to a ball $B(0,C)$, and may restrict $\omega$ to make an angle of $\ll 1$ with the vertical.  We now perform routine discretization.

\begin{definition}\label{point-def}
Define a \emph{(discretized) point} to be an element of the lattice $P := N^{-1} \Z^n \cap B(0,C)$.
Define a \emph{(discretized) line} to be the set of all points in $P$ which are within $C N^{-1}$  of a line, whose direction makes an angle of $\ll 1$ with the vertical.  Thus lines have cardinality $\sim N$, and their directions are defined up to uncertainty $O(N^{-1})$.  For any $N^{-1} \leq r \leq 1$, we say that a collection $\T$ of lines are \emph{$r$-separated} if their directions are $r$-separated.  If $r = N^{-1}$, we shall simply say that the collection $\T$ is \emph{separated}.
\end{definition}

Note that a separated collection of lines must have cardinality at most $\lesssim N^{n-1}$.  Unfortunately two points can determine more than one line; indeed we have 
\be{euclid}
\#( \{ T \in \T: x_1, x_2 \in T \} ) \lessapprox |x_1-x_2|^{1-n}
\end{equation}
for any separated collection $\T$ of lines.  This loss of $|x_1-x_2|^{1-n}$ motivates the two-ends reduction \eqref{two-ends} below.

By standard discretization arguments, the desired Kakeya estimate will follow from the discretized version
\be{strong}
(\sum_{T \in \T} (\sum_{x \in T} f(x))^{(n-1)p'})^{1/(n-1)p'} \lessapprox 
N^{1/p'} \| f \|_{l^p(P)}
\end{equation}
for all finite collections $\T$ of lines, no two of which are essentially parallel.  Note that we must have $\#(\T) \lessapprox N^{n-1}$ for such collections.

\begin{definition}\label{yt-def}
Let $\T$ be a separated collection of lines.  We define a \emph{shading} of $\T$ to be a map $Y: T \mapsto Y(T)$ on $\T$ such that $Y(T)$ is a subset of $T$ for all $T \in \T$.  
Let $0 < \lambda < 1$.  We say that the shading $Y$ has \emph{density $\lambda$ on $\T$} if
$$ \#(Y(T)) \approx \lambda \#(\T) \approx \lambda N \hbox{ for all } T \in \T.$$
If we replace $\approx$ by $\gtrapprox$, then we say that $Y$ has \emph{density at least $\lambda$ on $\T$}, etc.  We define the \emph{counting function} $\mu_{Y,\T}$ of $Y$ by
$\mu_{Y, \T}(x) := \sum_{T \in \T} \chi_{Y(T)}(x)$, and the \emph{mass} $\mass_{\T}(Y)$ to be the number
$\mass_\T(Y) := \sum_{T \in \T} \#(Y(T)) = \| \mu_{Y,\T}\|_{l^1}$.
If we have $\mass_{\T}(Y) \approx N \#(\T)$, we say that the shading is \emph{saturated in $\T$}.  
\end{definition}

To prove \eqref{strong}, it suffices by the usual restricted weak-type reduction (see \cite{wolff:kakeya}) to prove
\be{rwt-0}
N\lambda \#(\T)^{1/(n-1)p'} \lessapprox N^{1/p'} \#(\bigcup Y(\T))^{1/p}
\end{equation}
for all separated collections $\T$ of lines, all $1/N < \lambda \leq 1$, and all shadings $Y$ with density $\lambda$ on $\T$. We can rewrite \eqref{rwt-0} using $p = (4n+3)/7$ as
\be{rwt}
\#(\bigcup Y(\T)) \gtrsim \lambda^{(4n+3)/7} N \#(\T)^{\frac{4}{7}}.
\end{equation}

Fix $\T$, $\lambda$, $Y$.  We now apply another reduction from \cite{wolff:kakeya}, namely the \emph{two-ends reduction}.  This asserts that we may assume the condition
\be{two-ends}
\#(Y(T) \cap B(x,r)) \lessapprox \lambda r^\sigma N
\end{equation}
for all $0 < r \ll 1$, $T \in \T$ and some constant $\sigma > 0$ depending on $p$, $n$. See \cite{wolff:kakeya} for further details.
As we shall see, once one assumes the two-ends condition \eqref{two-ends} then one can improve the power of $\lambda$ in \eqref{rwt} substantially (the same phenomenon also occurs in \cite{wolff:kakeya}). 

Let $E$ denote the set $E := \bigcup Y(\T)$.  By dyadic pigeonholing we can find a subset $E' \subset E$ such that $\mu_{Y,\T} \sim \mass_{\T}(Y)/\#(E')$ on $E'$.
Fix $E'$.  Recall from Section \ref{notation-sec} that the map $\gamma_n: E' \to N^{-1} \Z \cap B(0,C)$ is defined by
$\gamma_n(x_1, \ldots, x_n) := x_n$.
The fibers of $\gamma_n$ are thus slices of $E'$.
By dyadic pigeonholing the size of the slices, there thus exists a set $S \subset N^{-1} \Z \cap B(0,C)$ of size 
\be{s-size}
|S| \approx 2^{-k} N
\end{equation}
for some $k \geq 0$ such that
\be{slice}
\#(E' \cap \gamma_n^{-1}(t)) \approx 2^k \#(E') / N
\end{equation}
for all $t \in S$.  As we shall see, the worst case shall be when $k = O(1)$.

Fix $k$, $S$.  Let $Y'$ denote the shading $Y'(T) := Y(T) \cap E' \cap \gamma_n^{-1}(S)$.  Integrating $\mu_{Y,\T}$ on $E' \cap \gamma_n^{-1}(S)$ using \eqref{s-size}, \eqref{slice} we see that
$\mass_{\T}(Y') \approx \mass_{\T}(Y) \approx \lambda N \#(\T)$.
One can then find a refinement $\T'$ of $\T$ such that 
$Y'$ has density $\lambda$ on  $\T'$.

We first use a basic ``two-slices'' argument (equivalent to Bourgain's ``bush'' argument) to dispose of a relatively easy case when $\lambda$ is small.  For each $T \in \T'$ and $t_1 \in S$ the two-ends condition gives
$$ \#(\{ t_2 \in \gamma_n(Y'(T)): |t_2 - t_1| \approx 1 \}) \gtrapprox \lambda N.$$
Summing over all $t_1$ and then over all $T'$, we obtain
\be{easy-two-ends}
\#(\{ (T,t_1, t_2) \in \T' \times S^2: t_1, t_2 \in \gamma_n(Y'(T)); |t_1 - t_2| \approx 1 \}) \gtrapprox \lambda^2 N^2 \#(\T).
\end{equation}
We may therefore find $t_1, t_2 \in S$ such that $|t_1 - t_2| \approx 1$ and
$$ \#(\{ T \in \T': t_1, t_2 \in \gamma_n(Y'(T)) \}) \gtrapprox 2^{2k} \lambda^2 \#(\T).$$
On the other hand, by \eqref{euclid} for any fixed $x_1 \in \gamma_n^{-1}$, $x_2 \in \gamma_n^{-1}(t_2)$, there are $\lessapprox 1$ lines $T$ which contain both $x_1$ and $x_2$.  Thus
$$ \#(\{ T \in \T': t_1, t_2 \in \gamma_n(Y'(T)) \}) \lessapprox \#(\gamma_n^{-1}(t_1)) \#(\gamma_n^{-1}(t_2))\lessapprox 2^{2k} \#(E)^2/N^2.$$
Combining these bounds one obtains $\#(E) \gtrapprox \lambda N \#(\T)^{1/2}$.
Since $p = (4n+3)/7$, $n > 8$, and $\#(\T) \lesssim N^{n-1}$, this bound will imply \eqref{rwt} when $\lambda \lessapprox N^{-1/8}$. 

\begin{figure}[htbp] \centering
 \ \psfig{figure=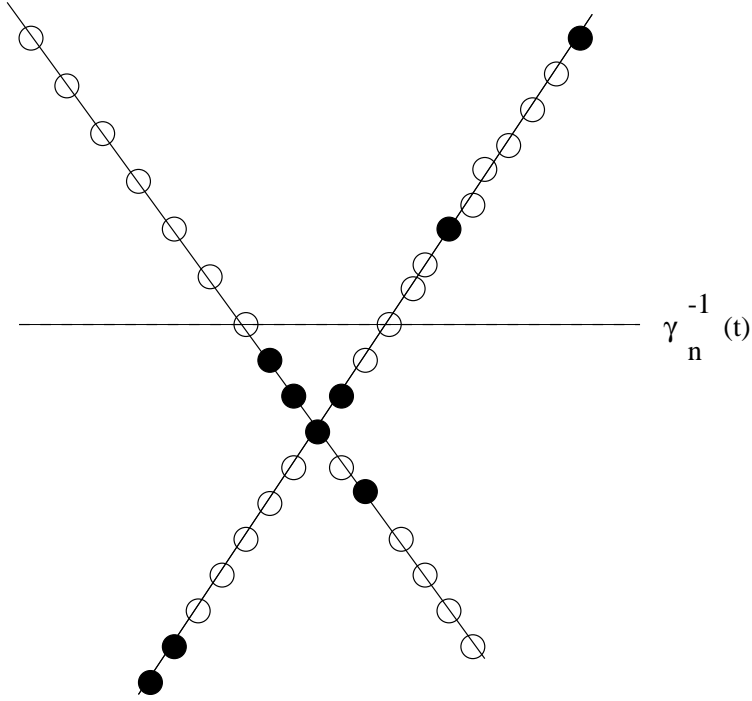}
 \caption{Two separated lines, and a shading of the lines (indicated by 
filled points as opposed to empty points).  Note that the NW-SE line fails to
satisfy the two-ends condition, while the NE-SW line obeys \eqref{four-ends}.
Note that if the density $\lambda$ is small, then many slices $\gamma_n^{-1}(t)$ will only have a small intersection with $E$.  This is the main difficulty
in applying the slices method to the maximal function problem.
}
 \end{figure}

To handle the remaining case 
$N^{-1/8} \lessapprox \lambda \lessapprox 1$ we 
use a ``six-slices'' argument. Certain ratios of these slices should be slopes of the type used in Theorem \ref{conviviality}; the flexibility we have in choosing these slopes is what allows us to get a sufficiently good power of $\lambda$ in our estimates.

By repeating the derivation of \eqref{easy-two-ends} we have
\be{four-ends}
\#(\{ (t_1,t_2,t_3,t_4) \in \gamma_n(Y'(T))^4: |t_i - t_j| \approx 1 \hbox{ for all } 1 \leq i < j \leq 4 \}) \gtrapprox \lambda^4 N^4
\end{equation}
for all $T \in \T'$.  On the other hand, $Y'$ has density $\lambda$ on $\T$.  If we then define the sets
$$ Q_{t_1,t_2}(T) := \{(t_3,t_4) \in \gamma_n(Y'(T))^2: |t_i - t_j| \approx 1 \hbox{ for all } 1 \leq i < j \leq 4 \}$$
and
$$ P(T) := \{ (t_1, t_2) \in \gamma_n(Y'(T))^2: \#(Q_{t_1,t_2}(T)) 
\approx \lambda^2 N^2 \} $$
we see that
$\#(P(T)) \approx \lambda^2 N^2$
if the implicit constants are chosen appropriately.  Summing this over all $T$, we obtain
$$ \sum_{t_1, t_2 \in S} \#(\{T \in \T': (t_1,t_2) \in P(T)\})
\approx \lambda^2 N^2 \#(\T);$$
by \eqref{s-size}, we can therefore find $t_1, t_2 \in S$ such that
\be{tpppp-card}
\#(\T'') \gtrapprox 2^{2k} \lambda^2 \#(\T),
\end{equation}
where $\T'' := \{T \in \T': (t_1,t_2) \in P(T)\}$.

Fix $t_1, t_2$, and define 
$S' := \{t \in S: |t-t_1|, |t-t_2| \approx 1\}$.
For any $t \in S'$, we define the associated slope $r(t)$ by $r(t) := (t-t_1)/(t_2-t)$, thus $|r(t)| \approx 1$. By construction, we have $Q_{t_1,t_2}(T) \subseteq S' \times S'$ and
$\#(Q_{t_1,t_2}(T)) \gtrapprox \lambda^2 N^2$ for all $T \in \T''$. 
Define the function $s$ on $S' \times S'$ by
$s(t,t^{\prime}):=r(t) +{r(t) \over r(t^{\prime})} + O(1/N)$,
where the $O(1/N)$ term is such that $s$ lands in $N^{-1}\Z \cap B(0,C)$.
From \eqref{cauchy} we have
$$ \#(\{ (t_3, t_4, t_5, t_6) \in Q_{t_1,t_2}(T)^2:  
(t_3, t_4) \sim_s (t_5, t_6) \}) \gtrapprox \lambda^4 N^3.$$
We now refine this to
$$ \#(\{ (t_3, t_4, t_5, t_6) \in Q_{t_1,t_2}(T)^2:  
(t_3, t_4) \sim_s (t_5, t_6); |t_3 - t_5| \gtrapprox \lambda^2 \}) \gtrapprox \lambda^4 N^3.$$
To see this, we count the exceptional set when $|t_3 - t_5| \ll \lambda^2$.
From the construction of $\T'$ there are $\lessapprox \lambda N$ choices of $t_3$, and similarly for $t_4$.  Fixing $t_3$, there are then $\ll \lambda^2 N$ choices of $t_5$.  Finally, for fixed $t_3$, $t_4$, $t_5$ there are only $\lessapprox 1$ choices for $t_6$, and so the exceptional set is suitably small.

By dyadic pigeonholing twice, we may thus find we may fix a $\lambda^2 \lesssim d \lesssim 1$ and a refinement $\T'''$ of $\T''$ such that
$$ \#(\{ (t_3, t_4, t_5, t_6) \in Q_{t_1,t_2}(T)^2:  
(t_3, t_4) \sim_s (t_5, t_6); |t_3 - t_5| \sim d \}) \gtrapprox \lambda^4 N^3$$
for all $T \in \T'''$.  Summing over $T$ using \eqref{tpppp-card}, we obtain
\bas \#(\{ (T, t_3, t_4, t_5, t_6) \in &\T''' \times (S')^4: (t_3, t_4), (t_5, t_6) \in Q_{t_1,t_2}(T); 
\\
(t_3, t_4) \sim_s (t_5, t_6); |t_3 - t_5| \sim d \})
& \gtrapprox \lambda^6 N^3
2^{2k}  \#(\T).
\end{align*}
On the other hand, the set
$$ \Delta := \{ (t_3, t_4, t_5, t_6) \in (S')^4: 
(t_3, t_4) \sim_s (t_5, t_6); |t_3 - t_5| \sim d \}$$
has cardinality $\lessapprox (2^{-k} N)^2 (dN)$ by \eqref{s-size} and the same counting argument used to prove \eqref{easy-two-ends}; note that $dN \gtrsim \lambda^2 N \gg 1$.  By the pigeonhole principle we may therefore find $(t_3, t_4, t_5, t_6) \in \Delta$ such that 
\be{t5-card}
\#(\T'''')
\gtrapprox \lambda^6 2^{4k} \#(\T) / d,
\end{equation}
where
$\T'''' := \{ T \in \T''': (t_3, t_4), (t_5, t_6) \in Q_{t_1,t_2}(T) \}$.

Fix $t_3, t_4, t_5, t_6, \T''''$.
We now run the argument from the proof of Theorem \ref{conviviality}, with 
$$G := \{ (a,b) \in \gamma_n^{-1}(t_1) \times \gamma_n^{-1}(t_2): a,b \in T \hbox{ for some } T \in \T'''' \},
$$
$(r_0, r_1, r_2, r'_1, r'_2, r_\infty) := (0, r(t_3), r(t_5), r(t_4), r(t_6), \infty),$
and $s := s(r_3,r_4)=s(r_5,r_6)$.  The reader may easily verify that
$\frac{s}{r_i} - \frac{1}{r'_i} = 1 + O(1/N)$ for $i=1,2,$
so that $r'_i$ is essentially the dual of $r_i$ in the sense of Section \ref{basic-sec}.  Also observe that $|r_1 - r_2| \approx d$, and $|r_1|, |r'_1|, |r_2|, |r'_2| \approx 1$.

Each line $T$ contributes $\approx 1$ elements to $G$, thus by \eqref{t5-card} we have
\be{g-card}
\#(G) \gtrapprox \lambda^6 2^{4k} \#(\T) / d.
\end{equation}
By the same token we see that the map $\pi_{-1}$ defined in Section \ref{basic-sec} is essentially one-to-one, up to a multiplicity of $\approx 1$. 

Also, by construction we have that $\pi_{r_1}(G)$ is essentially contained in a linear transformation of $\gamma_n^{-1}(t_3) \cap \bigcup \T'''$.  By construction of $\T'''$, this set is contained in $\gamma_n^{-1}(t_3) \cap E'$.  By \eqref{slice} we thus have
$ \# (\pi_{r_1}(G)) \lessapprox 2^k \#(E) / N$.
Similarly with $r_1$ replaced by $0$, $\infty$, $r'_1$, $r_2$, and $r'_2$.

We are almost ready to apply Theorem \ref{conviviality}.  However there is one remaining snag, namely that \eqref{co-ord} breaks down, or in other words knowledge of $\pi_r(g)$ and $\pi_{\tilde r}(g)$ do not necessarily determine $g$.  When $|r-\tilde r| \approx 1$ this is not an issue since one still has \eqref{co-ord} holding (modulo a multiplicity of $\approx 1$).  If one then inspects the proof of Theorem \ref{conviviality}, we see that all the arguments continue to work except for the stage where one asserts that $V$ is parameterized by $\pi_{r_1 \otimes r'_1}$ and $\pi_{r_2 \otimes r'_2}$.  Instead, for fixed values of $\pi_{r_1 \otimes r'_1}(v)$ and $\pi_{r_2 \otimes r'_2}$, the separation properties of the $r_i$ allow there to be as many as $\lessapprox d^{1-n}$ possible values of $v$ (because \eqref{euclid} gives this many choices for the first component $\gamma_1(v)$, and the second component $\gamma_2(v)$ is then essentially fixed thanks to the separation between $r'_1$, $r'_2$ and $0$).  By inspection of the proof of Theorem \ref{conviviality}, we see that this loss of one-to-oneness eventually leads to a loss of $d^{(1-n)/4}$ in the final upper bound for $\#(G)$.  We thus have
$\#(G) \lessapprox d^{(1-n)/4} (\frac{2^k \#(E)}{N})^{7/4}$.
Applying \eqref{g-card} and using $k \geq 0$, $d \gtrapprox \lambda^2$ we obtain  
$\#(E) \gtrapprox N \lambda^{\frac{2n+14}{7}} \#(\T)^{\frac{4}{7}}$.  Since $n > 8$, \eqref{rwt} follows.

\section{Hausdorff}\label{hausdorff-sec}

The purpose of this section is to prove an estimate on the Hausdorff dimension
of Kakeya sets. What we do here differs from the rest of the paper in that
we do not divide the set into slices. We essentially combine the
``iterated arithmetic techniques" of section 2 with the
``hairbrush" based ideas of \cite{wolff:kakeya}.

As is well known (see \cite{borg:kakeya}), the bound \eqref{rwt-0} on the maximal function at exponent $p$ implies that Besicovitch sets in $\R^n$ have dimension at least $p$.  Indeed, only needs to prove \eqref{rwt-0} in the special case $\lambda \approx 1$.  Heuristically, an application of $SD(1 + \frac{\sqrt{2}}{2})$ should obtain a dimension bound of $(2-\sqrt{2})(n-1)+1$.  We can improve this slightly to $(2-\sqrt{2})(n-4) + 3$ by eliminating the use of slices and work directly with the Besicovitch set, exploiting an old argument of Cordoba as in \cite{wolff:kakeya} to extract an additional gain.

As before, $N \gg 1$ is a large integer.  We shall also be working in some large ambient dimension $\R^M$ (not necessarily $\R^n$); the implicit constants may depend on $M$.  We adopt the notations of (discretized) points and lines from the previous section, with the following additional notations.
We say that two lines are \emph{essentially parallel} if their directions are within $O(1/N)$ of each other.
We say that a collection of points are \emph{essentially collinear} if they lie inside a common line.  
We say that a collection of points and lines are 
\emph{essentially coplanar} if they lie inside a $CN^{-1}$ neighbourhood of 
a two-dimensional plane. 

We define a \emph{$n-1$-collection of lines} to be any separated collection $\T$ of lines such that 
\be{l-angle}
\#(\{ T \in \T: \angle(T, \omega) \lesssim \theta \}) \lessapprox (N \theta)^{n-1}
\end{equation}
for all directions $\omega$ and angles $1/N \leq \theta \leq 1$, where $\angle(T,\omega)$ denotes the angle between $T$ and $\omega$ (defined up to an uncertainty of $O(1/N)$).  If in addition
 we have 
$\#(\T) \approx N^{n-1}$
then we say that $\T$ is \emph{saturated}.  

\begin{definition}\label{kakeya-def}  Let $n, d > 0$ be \emph{real} numbers.
We say that we have the \emph{Kakeya estimate} $K(n,d)$ if one has
$\#(\bigcup Y(\T)) \gtrapprox N^{d}$
whenever $\T$ is a saturated $n-1$-collection of lines in some ambient dimension $M$, and $Y$ is a saturated shading on $\T$.
\end{definition}

The statement $K(n,d)$ implies 
 that Kakeya sets in ${{\Bbb R}^n}$ have Hausdorff dimension at least $d$.
Also, $K(n,d)$ automatically implies the generalization
\be{kakeya-general}
\#(\bigcup Y(\T)) \gtrapprox \frac{\#(\T)}{N^{n-1}} N^d
\end{equation}
if the assumption that $\T$ is saturated is omitted.  This is by the usual factorization  argument, joining together random rotations of $L$ to create a saturated set (deleting the clusters of directions for which \eqref{l-angle} fails; note that these regions are generically small when $\#(\T) \ll N^{n-1}$). 

The purpose of this section is to show the functional relationship

\begin{theorem}\label{iter}  Let $0 < d < n$.
The $K(n,d)$ and $K(d+1, d')$ imply $K(n, \frac{2n + 1 + d'}{4})$.
\end{theorem}

Assume for the moment that Theorem \ref{iter} holds.  Then any statement of the 
form
$$
K(n, a(n-4) + 3 - b) \hbox{ for all } n \in \R^+
$$
for some constants $a, b > 0$ immediately implies the variant
$$ K(n, \frac{a^2+2}{4}(n-4) + 3 - b(\frac{a+1}{2})) \hbox{ for all } n \in \R^+.$$
(apply Theorem \ref{iter} with $d := a(n-4) + 3-b$ and $d' := a(d-3) + 3-b$). By iteration we thus have $K(n, (2-\sqrt{2})(n-4) + 3 + \eps)$ for all $\eps > 0$, which gives the desired Hausdorff bound.

\begin{proof}
Fix $n$, $d$, $d'$.  By raising $d$ if necessary we may assume that $K(n,d+\eps)$ fails.  It will then suffice to prove that $(2n+1+d')/4 \leq d+C\eps$, since the claim follows by letting $\eps \to 0$.

Since $K(n,d)$ holds and $K(n,d+\eps)$ fails, we may find an ambient dimension $M$, a saturated $n-1$-collection $\T$ of lines, and a saturated shading $Y$ such that
\be{e-card}
\#(E) \approx N^d,
\end{equation}
where $E := \bigcup Y(\T)$.
It then suffices to show that
\be{e-targ}
 \#(E) \gtrapprox N^{(2n+1+d')/4}.
\end{equation}

We now perform a reduction similar to the two-ends reduction, which ensures that $E$ has dimension $d$ in an appropriate sense.
Let $C$ be a large constant to be chosen later, and let $1/N < r < 1$ be a dyadic radius.  Call a ball $B(x,r)$ \emph{heavy} if one has
\be{heavy-def}
\#(E \cap B(x,r)) \geq N^{C\eps} (Nr)^{d}.
\end{equation}
Let $E_r$ denote the set of all points in $E$ which lie in at least one heavy ball $B(x,r)$.  We claim that
\be{sparse-clump}
\sum_{T \in \T} \#(Y(T) \cap E_r) \leq N^{-\eps} \sum_{T \in \T} \#(Y(T))
\end{equation}
if $C$ is sufficiently large.

To see this, suppose for contradiction that the above estimate failed.  Then the shading $Y_r(T) := Y(T) \cap E_r$ is saturated.  We may then find a refinement $\T_r$ of $\T$ such that $Y'$ has density $1$ on $\T_r$.

The set $E_r$ is contained in a union of heavy $r$-balls.  By applying the Kakeya hypothesis at eccentricity $1/r$, applied to an appropriate blurring of the collection $\T_r$ and the shading $Y_r$, we thus see that the number of heavy $r$-balls needed to cover $E_r$ is $\gtrapprox r^{-d}$.  From \eqref{heavy-def} we thus have
$ \#(E_r) \gtrapprox N^{C\eps} (rN)^d r^{-d},$
contradicting \eqref{e-card} if $C$ is sufficiently large.   This proves \eqref{sparse-clump}.  Henceforth the implicit constants may depend on $C$.

Set
$E' := E \backslash \bigcup_{1/N < r < 1} E_r$
and define the shading $Y'(T) := Y(T) \cap E'$ on $\T$.  Since $Y$ is saturated on $\T$, we see from \eqref{sparse-clump} that $Y'$ is also saturated on $\T$. 

Since $Y'$ is saturated on $\T$, we may find a refinement $\T'$ of $\T$ such that $Y'$ has density 1 on $\T'$, so in particular
$$ \# (\{ (x,T) \in E' \times \T': x \in Y'(T) \}) = \mass_{\T'}(Y') \approx N \#(\T') \approx N^{n}.$$
We now apply the ``bilinear reduction'' \cite{tvv:bilinear}. From \eqref{cauchy} and pigeonholing we have 
$$ \# (\{ (x, T_1,T_2) \in E' \times (\T')^2: x \in Y'(T_1) \cap Y'(T_2);
\angle(T_1,T_2) \sim \theta \}) \gtrapprox N^{2n-d}$$
for some $N^{-1} \leq \theta \leq 1$;
the diagonal contribution $l_1 = l_2$ can be discarded since $d < n$.  By another  pigeonholing we can find a direction $\omega$ such that
$$ \#(\{ (x, T_1,T_2) \in E' \times (\T')^2: x \in Y'(T_1) \cap Y'(T_2);
\angle(T_1,T_2) \sim \theta; \angle(T_1,\omega), \angle(T_2,\omega) \lesssim \theta \})\gtrapprox N^{2n-d} \theta^{n-1}.$$
As in \cite{tvv:bilinear} we can then rescale $\theta$ to equal 1 to obtain
\be{a-card}
\#(A) \gtrapprox N^{2n-d}.
\end{equation}
where $A := \{ (x, T_1,T_2) \in E' \times (\T')^2: x \in Y'(T_1) \cap Y'(T_2);
\angle(T_1,T_2) \sim 1 \}$.
We refer to elements $\alpha \in A$ as \emph{angles}, and write $x(\alpha), T_1(\alpha), T_2(\alpha)$ for $x$, $T_1$, $T_2$. 

\begin{definition}\label{pivot-def}
Let $\alpha$ be an angle.  We say that a point $z \in P$ is a \emph{pivot} for 
$\alpha$ if

\begin{itemize}
\item $z$ is distance $\approx 1$ from both $T_1(\alpha)$ and $T_2(\alpha)$, and is 
essentially coplanar with $T_1(\alpha)$, $T_2(\alpha)$.
\item There exists a point $i_{\alpha,z} \in Y'(T_1(\alpha))$ such that $c-i_{\alpha,z}$ 
is essentially parallel to $T_2(\alpha)$.
\item We have
\be{many-hit}
\#(\{ (y_1, y_2) \in Y'(T_1(\alpha)) \times Y'(T_2(\alpha)): y_1, z, y_2 \hbox{ are 
essentially collinear} \}) \approx N.
\end{equation}
\end{itemize}
\end{definition}

\begin{figure}[htbp] \centering
 \ \psfig{figure=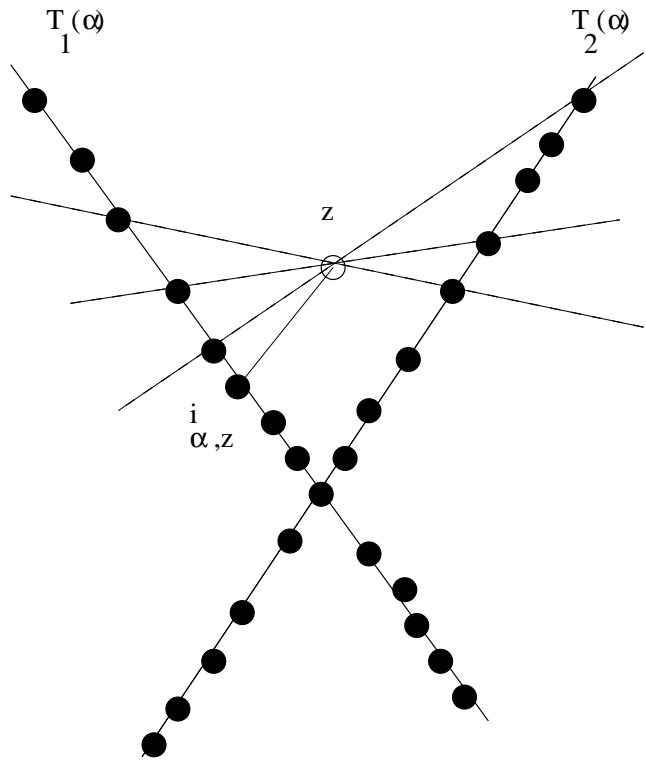}
 \caption{An angle and one of its pivots.  For the Hausdorff problem the shadings have density $\approx 1$.}
\end{figure}

Angles and pivots are the analogue of vertical line segments and values of $\nu$ in Section \ref{basic-sec}.
Let 
$\Omega := \{ (\alpha, z) \in A \times P: z \hbox{ is a pivot for } \alpha \}$
denote the set of angle-pivot pairs.
The following lemma asserts that a large fraction of the points coplanar to $T_1(\alpha)$  and $T_2(\alpha)$ are indeed pivots for $\alpha$.

\begin{lemma}\label{pivot-lemma}
For every $\alpha \in A$ we have
$\#(\{ z \in P: (\alpha,z) \in \Omega \}) \approx N^2.$
\end{lemma}

\begin{proof}
The upper bound is trivial, so it suffices to prove the lower bound.
Fix $\alpha$, and let $X$ denote the set of triples
$$ X := \{ (i, y_1, y_2) \in Y'(T_1(\alpha))^2 \times Y'(T_2(\alpha)):
 |i-y_1|, |i-x(\alpha)|, |y_1-x(\alpha)|, |y_2-x(\alpha)| \approx 1 \}.$$
Since $\#(Y'(T_1(\alpha))^2 \times Y'(T_2(\alpha)))\approx N^3$
we see that
$\#(X) \approx N^3$ (if the constants are chosen correctly).

For each triple $\tau = (i, y_1, y_2) \in X$, we can find a point $\nu(\tau)$ which is
 distance $\approx 1$ from $T_1(\alpha)$ and $T_2(\alpha)$, is essentially coplanar
 with $T_1(\alpha), T_2(\alpha)$, is essentially collinear with $y_1$, $y_2$, and is  such that $\nu(\tau) - i$ is essentially parallel to $T_2(\alpha)$.  Let us fix this map $\nu: X \to P$.  

Since the $\nu(\tau)$ are essentially coplanar to $T_1(\alpha)$ and $T_2(\alpha)$ we see that $\#(\nu(X)) \lessapprox N^2$.  On the other hand, we see from elementary geometry that
$\#(\nu^{-1}(z)) \lessapprox N$ for all $z \in \P$, and thus 
$$ \#(\{ z \in \nu(X): \#(\nu^{-1}(z)) \approx N \}) \approx N^2.$$
Since all the elements $z$ in this set are pivots for $\alpha$, the claim follows.
\end{proof}

From the preceding Lemma and \eqref{a-card} we thus have
\be{om-card}
\#(\Omega) \approx N^2 \#(A) \gtrsim N^{2n - d + 2}.
\end{equation}

Let $f: \Omega \to P \times P$ denote the map $f(\alpha,z) := (z, i_{\alpha,z})$.

\begin{lemma}\label{cordoba}
We have $\#(f(\Omega)) \approx N^2 \#(A) \approx \#(\Omega)$.
\end{lemma}

\begin{proof}
This will be a variant of Cordoba's argument.  We first observe that $f(\alpha,z)$ 
essentially determines $T_2(\alpha)$ (up to a multiplicity of $\approx 1$), since $T_2(\alpha)$ is parallel to $z - i_{\alpha,z}$, which has 
magnitude $\approx 1$, and the $T_2(\alpha)$ are direction-separated.  (This is the analogue of $\nu$ being determined by $\pi_{r_\infty} \circ \gamma_1$ in Section \ref{basic-sec}).  Thus we may find a refinement $\Omega'$ of $\Omega$
such that $T_2(\alpha)$ is determined by $f(\alpha, z)$ on $\Omega'$.

By \eqref{cauchy} and a pigeonholing it suffices to show that
$$ \#(\{ ((\alpha, z), (\alpha', z')) \in \Omega' \times \Omega': (\alpha,z) \sim_f (\alpha',z'); \angle(T_1(\alpha), T_1(\alpha')) \sim \theta \}) \lessapprox N^2 \#(A)$$
for each $1/N \leq \theta \leq 1$.

Fix $\theta$, and consider an element $\alpha \in A$.  If $((\alpha,z), (\alpha', z))$ contributes to the above set, then $z=z'$ and $T_1(\alpha)$,
 $T_1(\alpha')$, and $T_2(\alpha) = T'_2(\alpha)$ are essentially coplanar, so the number of possible $T_1(\alpha')$ is $\lessapprox N\theta$ since $L$ is separated.  For each such $T_1(\alpha')$, the 
number of $i_{\alpha,z}$ which can contribute is 
$\leq \#(T_1(\alpha) \cap T_1(\alpha')) \lessapprox 1/\theta$.
 Since $z - i_{\alpha,z}$ is essentially parallel to $T_2(\alpha)$, 
the number of $z$ which can contribute is $\lessapprox N/\theta$.  The claim follows. 
\end{proof}

We use this lemma to fix a refinement $\Omega'$ of $\Omega$ such that $\Omega'$ is parameterized by $f$.  We now lift our Besicovitch set to the space $\R^{M+1} := \{ (x,t): x \in \R^M, 
t \in \R \}$ (this is the analogue of the iteration argument in Section \ref{basic-sec}).  If $(\alpha,z) \in \Omega'$, we define $\Tlift(\alpha,z) \subseteq 
\R^{M+1}$ to be a (discretized)  line which contains the points
$(x(\alpha),0)$ and $(x(\alpha) + (z - i_{\alpha,z}), 1)$, and whose $t$ variables ranges over the region $|t| \lessapprox 1$. For each $z \in P$, define $\Tblift(z) := \{ \Tlift(\alpha,z): (\alpha,z) \in \Omega' \}$.  We also
 define $\Ylift(\Tlift(\alpha,z))$ to be those elements $(y_1,t)$ of
 $\Tlift(\alpha,z)$ such that $y_1$ appears in the left-hand side of \eqref{many-hit}.  By construction, $\Ylift$ is a shading with density 1 on $\Tblift(z)$.
From elementary geometry we see that $|t - 1| \approx 1$ for all $(x,t) \in \Ylift(\Tlift(\alpha,z))$.

\begin{lemma}\label{card-bound} $\Tblift(z)$ is a $d$-collection of lines for each $z \in P$. 
\end{lemma}

\begin{proof}
Fix $z$.  Since the direction of $\Tlift(\alpha,z)$ is affinely determined by $i_{\alpha,z}$, it suffices to show 
$$
\#(\{ \alpha \in A: (\alpha, z) \in \Omega'; i_{\alpha,c} \in B(x, \theta) \})
\lessapprox (N \theta)^{d-1}
$$
for all balls $B(x,\theta)$.  On the other hand, from the construction of $E'$ we have
$\#(E' \cap B(x,r)) \lessapprox (Nr)^{d}$. Since $\Omega'$ is parameterized by $f$, the claim follows.
\end{proof}

From this lemma and the hypothesis $K(d+1,d')$ we can now obtain a multiplicity 
bound on $\Ylift$.  More precisely: 

\begin{lemma}\label{upper}  Fix $z \in P$.  Then there exists a set of points $\overline{E}(z) \subset \R^{M+1}$ such that the shading $\Ylift'(\Tlift(\alpha,z)) := \Ylift(\Tlift(\alpha,z)) \cap \overline{E}(z)$
is saturated on $\Tblift(z)$, and one has $\mu_{\Ylift', \Tblift(z)} \lesssim N^{d+1-d'}$ on $\overline{E}(z)$.
\end{lemma}

\begin{proof}
We set
$\overline{E}(z) := \{ (x,t) \in \R^{M+1}: \mu_{\Ylift', \Tblift(z)}(x,t) \lessapprox N^{C\eps} N^{d+1 - d'}\}$
where $C$ is to be chosen later. 
Now suppose for contradiction that $\Ylift'$ is not saturated.  Then the shading 
$\Ylift''(\Tlift(\alpha,z)) :=  \Ylift(\Tlift(\alpha,z)) \backslash \overline{E}(z)$
must be saturated on $\Tblift(z)$.  There must therefore be a refinement $\Tblift'(z)$ of $\Tblift(z)$ such that $\Ylift''$ has density 1 on $\Tblift'(z)$.  By the hypothesis $K(d+1,d')$ we thus have
$$ \#(\bigcup \Ylift''(\Tblift'(z))) 
\gtrapprox N^{-d} \#(\Tblift'(z)) N^{d'}
\approx N^{d'-d} \#(\Tblift(z)).$$
By construction, $\bigcup \Ylift''(\Tblift'(z))$ is outside $E$.  From the definition of $E$ we thus have 
$$ \| \sum_{\Tlift(\alpha,z) \in \Tblift(z)} \chi_{\Tlift(\alpha,z)} \|_1 \gtrapprox N^{C\eps} N^{d+1-d'} N^{d'-d} \#(\Tblift(z)).$$
On the other hand, the left-hand side is clearly $\lessapprox N \#(\Tblift(z))$.  Thus we have a contradiction if $C$ is sufficiently large. 
\end{proof}

Fix $z \in P$, and let $\Ylift'$ be as in the previous lemma.  The counting
function $\mu_{\Ylift', \Tblift(z)}$ 
has an $l^1$ norm of $\lessapprox N \#(\Tblift(z))$ and an $l^\infty$ norm of $\lessapprox N^{d+1-d'}$.  From H\"older and a summation $z$ we thus have the lower bound
$\sum_{z \in P} \| \mu_{\Ylift', \Tblift(z)} \|_{l^2}^2
\lessapprox N^{d+2 - d'} \#(\Omega')$.
To prove \eqref{e-targ}, it thus suffices by \eqref{om-card} to obtain the upper bound
\be{2-lower}
\sum_{z \in P} \| \mu_{\Ylift', \Tblift(z)} \|_{l^2}^2
\gtrapprox N^{1 - 2d} \#(\Omega')^2.
\end{equation}
We expand the left-hand side of \eqref{2-lower} as
$$
\sum_{x \in P} \sum_{t \in N^{-1} \Z} \sum_{(\alpha,z), (\alpha',z) \in \Omega'} 
\chi_{\Ylift'(\Tlift(\alpha,z))}(x,t) \chi_{\Ylift'(\Tlift(\alpha',z))}(x,t).$$

Suppose $(\alpha, z) \in \Omega'$, $x_1 \in Y'(T_1(\alpha))$, $x_2 \in Y'(T_2(\alpha))$ are 
such that $x_1$, $z$, and $x_2$ are essentially collinear.  Then there exists a 
$t \in N^{-1} \Z$ such that $(x_1,t) \in \Ylift'(\Tlift(\alpha,z))$.  In light of this, the
previous sum can be bounded from below by
$$
\sum_{(x_1, z, x_2) \in \Sigma}
\sum_{\alpha, \alpha' \in A} \chi_{\Omega'}(\alpha,z) \chi_{\Omega'}(\alpha',z)
\chi_{Y'(T_1(\alpha)) \cap Y'(T_1(\alpha'))}(x_1) \chi_{Y'(T_2(\alpha)) \cap
Y'(T_2(\alpha'))}(x_2),$$
where
$$ \Sigma := \{ (x_1, z, x_2) \in E \times P \times E: x_1, z, x_2 \hbox{ essentially 
collinear}; |x_1 - z|, |x_2 - z| \approx 1 \}.$$
The previous sum can be simplified to
$$
\sum_{(x_1, z, x_2) \in \Sigma}
 [\sum_{\alpha \in A} \chi_{\Omega'}(\alpha,z)
\chi_{Y'(T_1(\alpha))}(x_1) \chi_{Y'(T_2(\alpha))}(x_2)]^2.$$
By the Cauchy-Schwarz inequality, this can be bounded from below by
$$
[
\sum_{(x_1, z, x_2) \in \Sigma}
\sum_{\alpha \in A} \chi_{\Omega'}(\alpha,z)
\chi_{Y'(T_1(\alpha))}(x_1) \chi_{Y'(T_2(\alpha))}(x_2)]^2 / \#(\Sigma).
$$
From the definition of $\Sigma$ it is clear that
$\#(\Sigma) \approx N \#(E)^2.$
Also, for fixed $(\alpha, z) \in \Omega'$ we see from the definition of pivot that
$$ \sum_{x_1, x_2: (x_1, z, x_2) \in \Sigma} \chi_{Y'(T_1(\alpha))}(x_1) 
\chi_{Y'(T_2(\alpha))}(x_2) \approx N.$$
Thus the left-hand side of \eqref{2-lower} is bounded from below by
$N \#(E)^{-2} (\sum_{(\alpha, z) \in \Omega'} 1)^2,$
and the claim follows from \eqref{e-card}.
\end{proof}

\end{document}